\documentclass[12pt,a4paper]{article}

\usepackage{amsmath}
\usepackage{amssymb}
\usepackage{amsthm}

\newcommand{\CC}{\mathbb{C}} 
\newcommand{\RR}{\mathbb{R}} 
\newcommand{\ZZ}{\mathbb{Z}} 
\newcommand{\PP}{\mathbb{P}} 

\newcommand{\Hh}{\mathfrak{H}}

\newcommand{\M}{\mathcal{M}}

\newcommand{\T}{\mathcal{T}}
\newcommand{\K}{\mathcal{K}}

\newcommand{\I}{\mathcal{I}}

\newcommand{\tp}{{}^t}
\newcommand{\kuno}{\kappa}

\newcommand{\kenne}{\kappa}

\newtheorem{theorem}{Theorem}[section]
\newtheorem{proposition}[theorem]{Proposition}
\newtheorem{lemma}[theorem]{Lemma}

\newtheorem{definition}[theorem]{Definition}

\theoremstyle{remark}
\newtheorem{remark}[theorem]{Remark}

\numberwithin{equation}{section}

\DeclareMathOperator{\Sym}{Sym} \DeclareMathOperator{\im}{Im}
\DeclareMathOperator{\sgn}{sgn} \DeclareMathOperator{\End}{End}

\def\be{\begin{equation}}
\def\ee{\end{equation}}

\begin{document}

\title{\Large Vector-Valued Modular Forms from the Mumford Forms, Schottky-Igusa Form, Product of Thetanullwerte and the Amazing Klein Formula}

\author{Marco Matone \and Roberto Volpato}

\maketitle

\begin{abstract} Vector-valued Siegel
modular forms are the natural generalization of the classical
elliptic modular forms as seen by studying the cohomology of the
universal abelian variety. We show that for $g\geq4$, a new class of
vector-valued modular forms, defined on the Teichm\"uller space,
naturally appears from the Mumford forms, a question directly related
to the Schottky problem. In this framework we show that the
discriminant of the quadric associated to the complex curves of
genus $4$ is proportional to the square root of the products of
Thetanullwerte $\chi_{68}$, which is a proof of the recently
rediscovered Klein ``amazing formula''. Furthermore, it turns out
that the coefficients of such a quadric are derivatives of the
Schottky-Igusa form evaluated at the Jacobian locus, implying new
theta relations involving the latter, $\chi_{68}$ and the theta
series corresponding to the even unimodular lattices $E_8\oplus E_8$
and $D_{16}^+$. We also find, for $g=4$, a functional relation
between the singular component of the theta divisor and the Riemann
period matrix.
\end{abstract}

\newpage
\section{Introduction}

Here we shortly describe the results and introduce  some simple
(linear) algebraic facts that will be useful in the following and
describe our notation for theta functions.

For each fixed positive integers $g,n$, define  $$M_n:={g+n-1\choose
n}\ ,\; N_n:=(2n-1)(g-1)+\delta_{n1}\ ,\; K_n:=M_n-N_n \ ,$$ so that, for a
complex curve $C$ of genus $g\ge 2$, $M_n$ and $N_n$ are the
dimensions of $\Sym^n H^0(K_C)$ and $H^0(K_C^n)$, respectively. Let
$$\Hh_g:=\{Z\in M_g(\CC)\mid \tp Z=Z,\im Z>0\} \ ,$$ be the Siegel upper
half-space, i.e. the space of $g\times g$ complex symmetric matrices
with positive definite imaginary part, and define the usual action
of the symplectic group $\Gamma_g:={\rm Sp}(2g,\ZZ)$ on $\Hh_g$ by
$Z\mapsto (AZ+B)(CZ+D)^{-1}$, $\Big(\begin{matrix}A & B\\ C
&D\end{matrix}\Big)\in \Gamma_g$. Denote by $\M_g$ the moduli stack
of complex curves of genus $g$ and by $\mathcal{A}_g=\Gamma_g\backslash
\Hh_g$ the moduli stack of principally polarized abelian $g$-folds.
According to Torelli's Theorem the morphism ${\rm Jac}: \M_g \to
\mathcal{A}_g$, which on points takes the algebraic curves
to its Jacobian, is injective. The question of characterizing the
image of ${\rm Jac}$ is the Schottky problem.

Let $\mathcal{C}_g \stackrel{\pi}{\longrightarrow}\mathcal{M}_g$ be the
universal curve over  $\mathcal{M}_g$ and $L_n=R\pi_*(K^n_{\mathcal{C}_g/\mathcal{M}_g})$ the vector bundle on $\mathcal{M}_g$ of rank
$N_n$ with fiber $H^0(K_C^n)$ at the point of
$\mathcal{M}_g$ representing $C$. On $\mathcal{A}_g$ it is the conormal
bundle to the zero section in the universal abelian scheme; it pulls
back via ${\rm Jac}$ to $\M_g$. A vector-valued Siegel modular form
on $\mathcal{A}_g$ is a global section of some tensor
bundle of the Hodge vector bundle on $\mathcal{A}_g$, whose pull-back to $\M_g$ via the Torelli map corresponds to
$L_1$.
In spite
of their relevance, {\it e.g.} in number theory, they have been
studied essentially only in the case of genus two, where correspond to suitable
commutators of Siegel modular forms. In Section 2 we will introduce
a new set of vector-valued modular forms, defined on the
Teichm\"uller space, strictly related to the Mumford forms, which
are holomorphic global sections of the bundle
$${\Bbb F}_g=(\det {\Bbb E}_g)^{c_n}\otimes(\wedge^{N_n} \Sym^n{\Bbb E}_g)^* \ , $$
on $\M_g$, where ${\Bbb E}_g\equiv L_1$ and $c_n:=6n^2-6n+1$.  Their weight is
$c_n-{g+n-1\choose n-1}$.

Such vector-valued Teichm\"uller modular forms provide a deep
relation between Schottky's problem, quadrics describing $C$ and
$\theta$-functions. The construction fits into the general problem
of finding explicit formulae that relate theta nulls,
which can be
regarded as homogeneous coordinates on the moduli space of abelian varieties (with appropriate level structure), to projective
invariants of curves, which are homogeneous coordinates on $\M_g$. Such formulae are also very interesting for
physics, since they appear naturally in the computation of higher loop superstring
amplitudes \cite{Matone:2005vm,Matone:2008td,Matone:2010yv,Matone:2005bx,Matone:2012wy}.
As a first application of the vector-valued Teichm\"uller modular
forms, in Section 3 we consider the case $g=4$. In particular, we
will show that the discriminant of the quadric associated to $C$ is
proportional to the square root of the products of Thetanullwerte
$\chi_{68}$, which is a proof of the recently rediscovered Klein
``amazing formula''. The relation with the Schottky's problem is due
to the fact that the coefficients of such a quadric are derivatives
with respect to the period matrix of the Schottky-Igusa form. This
provides new theta relations involving the latter, $\chi_{68}$ and
the theta series corresponding to the even unimodular lattices
$E_8\oplus E_8$ and $D_{16}^+$. Finally, we will find a functional
relation between the singular component $\Theta_s$ of the theta
divisor and the Riemann period matrix and show that the
hyperelliptic locus is a zero of order $10$ for the Hessian of
$\theta(e,Z)$, $e\in\Theta_s$.

Given a basis $v_1,\ldots,v_g$ of a $g$-dimensional vector space
$V$, denote by $\tilde v_1^{(n)},\ldots,\tilde v_{M_n}^{(n)}$ the
basis of the symmetrized tensor product $\Sym^nV$ given by elements
of the form \be\label{symm}\frac{1}{n!}\sum_{\pi\in
S_n}v_{\pi(i_1)}\otimes\cdots\otimes v_{\pi(i_n)}\ ,\ee with $S_n$
the group of permutations of $n$ objects, taken with respect to an
arbitrary ordering.

\begin{proposition}\label{endo}Let $V\cong\CC^g$ be a
$g$-dimensional complex vector space, and fix $A\in GL(V)$; then,
the induced endomorphism on the $1$-dimensional space
$\wedge^{M_n}(\Sym^n V)\cong\CC$ is given by $\det A^{g+n-1\choose
n-1}\in GL_1(\CC)$. Explicitly, if $\tilde
w_i=\sum_{j=1}^gA_{ij}\tilde u_j$, $u\in V$, then
\be\label{wedsym}\tilde w^{(n)}_1\wedge\cdots\wedge \tilde
w^{(n)}_{M_n}=\det A^{\frac{n}{g}M_n}\,\tilde
u^{(n)}_1\wedge\cdots\wedge \tilde u^{(n)}_{M_n}\
.\ee\end{proposition}

Consider the theta function with characteristics
\be\label{thetaconc}\theta \left[^a_b\right]\left(z,Z\right):=
\sum_{k\in {\ZZ}^g}e^{\pi i \tp{(k+a)}Z(k+a)+ 2\pi i \tp{(k+a)}
(z+b)} \ ,\ee where $z\in \CC^g$, $Z\in \Hh_g$ and $a,b\in{\RR}^g$.
Let $\{\alpha_1,\ldots,\alpha_g,\beta_1,\ldots,\beta_g\}$ be a
symplectic basis of $H_1(C,\ZZ)$. Denote by $\{\omega_i\}_{1\le i\le
g}$ the basis of $H^0(K_C)$ satisfying the normalization condition
$\oint_{\alpha_i}\omega_j=\delta_{ij}$, and by
$\tau_{ij}:=\oint_{\beta_i}\omega_j\in\Hh_g$ the Riemann period
matrix, $i,j=1,\ldots,g$. The basis of $H_1(C,\ZZ)$ is determined up
to the
transformation \be\label{anchequesta}\begin{pmatrix}\alpha\\
\beta\end{pmatrix}\mapsto \begin{pmatrix}\tilde\alpha\\
\tilde\beta\end{pmatrix}=\begin{pmatrix}D & C \\ B& A\end{pmatrix}
\begin{pmatrix}\alpha\\ \beta\end{pmatrix}\ ,\qquad\qquad \gamma\equiv\begin{pmatrix}A
& B \\ C & D\end{pmatrix}\in \Gamma_g \ ,\ee which induces the
following transformation on the period matrix \be\label{modull}\tau
\mapsto \gamma\cdot\tau=(A\tau+B)(C\tau+D)^{-1}\ .\ee

For $n\in\ZZ$, denote by $J_n(C)$ the principal homogeneous space of
linear equivalence classes of divisors of degree $n$ on $C$. The
Jacobian $J(C):={\CC}^g/L_\tau$, $L_\tau:={\ZZ}^g +\tau {\ZZ}^g$, is
identified with $J_0(C)$: each point of $J_0(C)$ can be expressed as
$D_2-D_1$, with $D_1$ and $D_2$ effective divisors of the same
degree, which corresponds to  $\int_{D_1}^{D_2}\omega\in J(C)$.
Choose an arbitrary point $p_0\in C$ and let
$A(p):=(A_1(p),\ldots,A_g(p))$, $A_i(p):=\int_{p_0}^p\omega_i$,
$p\in C$, be the Abel-Jacobi map. It embeds $C$ into the Jacobian
$J_0(C)$ and generalizes to a map from the space of divisors of $C$
into $J_0(C)$ as $A(\sum_i n_i p_i):=\sum_in_iA(p_i)$, $p_i\in C$,
$n_i\in\ZZ$. By Jacobi Inversion Theorem the restriction of $A$ to
the space $C_g$ of divisors of degree $g$ on $C$ is a surjective map
onto $J_0(C)$. Consider the vector of Riemann constants
$K^p_i:=\frac{1}{2}+\frac{1}{2}\tau_{ii}-\sum_{j\neq
i}^g\oint_{\alpha_j}\omega_j\int_{p}^x\omega_i$, $i=1,\ldots,g$, for
all $p\in C$. For any $p\in C$ define the formal sum
$\Delta:=(g-1)p-K^p$ so that, for any divisor $\xi$ of degree $g-1$
in $C$, $\xi-\Delta$ is the point in $\CC^g$
 given by $\int_{(g-1)p}^\xi \omega+K^P$. Under the projection $\CC^g\to J_0(C)$, $\Delta$ becomes a distinguished point in
 $J_{g-1}(C)$ depending only on the homological class of the marking (recall that a marking for $C$ is given by fixing a canonical homotopy basis together with a basepoint $p_0\in C$,
 see e.g. \cite{FayMAM}).
 Furthemore,  $2\Delta =
K_C$. We refer to \cite{jfayy} and \cite{FayMAM} for further
details.

If $\delta',\delta''\in\{0,1/2\}^g$, then $\theta
\left[\delta\right]\left(z,\tau\right):=\theta
\left[^{\delta'}_{\delta''}\right]\left(z,\tau\right)$ has definite
parity in $z$ $\theta \left[\delta\right]
\left(-z,\tau\right)=e(\delta) \theta \left[\delta\right]
\left(z,\tau\right)$, where $e(\delta):=e^{4\pi i\!\tp{\delta'}
\delta''}$. There are $2^{2g}$ different characteristics of definite
parity. By Abel Theorem each one of such characteristics determines
the divisor class of a spin bundle $L_\delta\simeq K^{1\over2}_C$,
so that we may call them spin structures. There are $2^{g-1}(2^g+1)$
even and $2^{g-1}(2^g-1)$ odd spin structures.

We will also consider the prime form $E(z,w)$ and the multi-valued
$g/2$-differential $\sigma(z)$ on $C$ with empty divisor, satisfying
the property
$$
\sigma(z+\tp\alpha n+\tp\beta m)=\chi^{-g}e^{\pi i(g-1)\tp m \tau
m+2\pi i \tp m \K^z}\sigma(z)\ .
$$
Such conditions fix $\sigma(z)$ only up to a factor independent of
$z$; the precise definition, to which we will refer, can be given,
following \cite{FayMAM}, on the universal covering of $C$  (see also
\cite{jfayy}).

\begin{proposition}\label{thdettheta}For each integer $n$, let $\phi^n:=\{\phi_i^n\}_{1\le i\le N_n}$ be an arbitrary basis of $H^0(K_C^n)$.
Then \be\label{dettheta}\kuno[\phi^1]:=
{\det\phi_i^1(p_j)\sigma(y)\prod_1^gE(y,p_i)\over
\theta\bigl(\sum_{1}^gp_i-y-\Delta\bigr)\prod_1^g\sigma(p_i)
\prod_{i<j}^gE(p_i,p_j) }\ ,\ee for all $p_1,\ldots,p_g,y\in C$, and
\be\label{detthetaii}\kenne[\phi^n]:={\det \phi_i^{n}(p_j)\over
\theta\bigl(\sum_{1}^{N_n}
p_i-(2n-1)\Delta\bigr)\prod_{1}^{N_n}\sigma(p_i)^{2n-1}\prod_{i<j}^{N_n}
E(p_i,p_j)}\ ,\ee  for $n\geq2$, for all $p_1,\ldots,p_{N_n}\in C$,
depend only on the marking of $C$ and on $\{\phi_i^n\}_{1\le i\le
N_n}$.\end{proposition}

\vskip 3pt

\noindent {\sl Proof.} For each integer $n$, $\kenne[\phi^n]$ is a
meromorphic function on $C$ with empty divisor
\cite{jfayy,MatoneBB}. \hfill$\square$

\section{Vector-valued Teichm\"uller modular forms from the Mumford form}\label{secSiegel}

Let $\lambda_n:=\det L_n$ be the determinant line bundle. According
to Mumford \cite{Mumford}
$$
\lambda_n\cong\lambda_1^{\otimes c_n}\ ,
$$
where $c_n:=6n^2-6n+1$. The Mumford form $\mu_{g,n}$ is the unique,
up to a constant, holomorphic section of
$\lambda_n\otimes\lambda_1^{-\otimes c_n}$ nowhere vanishing on
$\mathcal{M}_g$.

Comparing $\mu_{g,2}$ with the Polyakov measure for the bosonic
string, Manin observed that $c_2=13$ in Mumford's formula coincides
with the half of the string critical dimension. In a seminal paper
\cite{BelavinCY} Belavin and Knizhnik proved that the Polyakov
measure coincides with $|\mu_{g,2}|^2$. More generally $-c_n$ is the
central charge of the chiral $b-c$ system of conformal weight $n$
\cite{BonoraCJ}.

Belavin and Knizhnik obtained $\mu_{g,2}$ from an expression for the
curvature form of the determinant of Laplace operators. As observed
in \cite{BostUY}, this is a particular case of the similar formula
for the determinants of Dirac operators on arbitrary compact
manifolds, due to Bismut and Freed \cite{BismutWR} (see also
\cite{AlvarezGaumeVM} and references therein). Such results lead to
expressions in terms of complex geometry of the canonical curve $C$
providing a link with the spectral invariants which appear using the
formulae for the Laplace operator determinants by Ray and Singer
\cite{RaySinger} leading to sums over lengths of closed geodesics by
means of the Selberg trace formula.

The expression of $\mu_{g,2}$ in terms of $\theta$-functions has
been  derived in the context of string theory by Beilinson and Manin
in \cite{BeilinsonZW}.

\vskip 3pt

\begin{theorem}\label{mumfordformm}Let $\{\phi^n_i\}_{1\le i\le N_n}$ be a basis of
$H^0(K_C^n)$, $n\geq2$. The Mumford form is, up to a universal
constant \be\label{diciamolo}
\mu_{g,n}={\kuno[\omega]^{(2n-1)^2}\over
\kenne[\phi^n]}{\phi^n_1\wedge\cdots\wedge\phi^n_{N_n}\over
(\omega_1\wedge\cdots\wedge\omega_g)^{c_n}}\ . \ee
\end{theorem}

\vskip 3pt

\noindent The expression of $\mu_{g,n}$ in terms of
$\theta$-functions,
given by Verlinde and Verlinde \cite{VerlindeKW} and Fay
\cite{FayMAM}, follows immediately by
 \eqref{dettheta} and \eqref{detthetaii}. Nevertheless, it remains the hard question of expressing
$\mu_{g,n}$ without using points on $C$. As we will see, there are
exceptions for $\mu_{2,2}$ and $\mu_{3,2}$. There is also a proof
for the long-standing conjecture for $\mu_{4,2}$, here we shows that
it follows immediately in the present context. A related issue concerning ${\kuno[\omega]^{(2n-1)^2}/\kenne[\phi^n]}$ is its dependence on the basis $\{\phi^n_i\}_{1\le i\le N_n}$: apparently, there is no natural choice for such a basis.
 In the case of
$\mu_{2,2}$ and $\mu_{3,2}$, since $K_2=0$ (i.e., $\Sym^2 H^0(K_C)$
and $H^0(K_C^2)$ have the same dimension), the natural choice is to
use $\{\omega_i\omega_j\}_{1\le i,j\le g}$ as basis of $H^0(K^2_C)$.
It turns out that such a choice reveals new interesting properties
just in the case when  $K_n>0$. Actually, the application of symmetric products of the $\omega_i$ 's in this case introduces free vectorial indices and therefore leads to vector-valued modular forms defined on the Teichm\"uller space (see \cite{vanderGeer} for a very nice account on vector-valued Siegel modular forms) which are strictly related to the investigations in \cite{MatoneBB,MatoneZZ,ShepherdBarronHC}. Remarkably, such a
structure will also lead to a strict connection between Mumford
forms, quadrics describing canonical curves, their discriminant and
the Schottky problem. In particular, for $g=4$ we will get some new
results connecting the above structures to the Schottky-Igusa form,
the theta series and the products of Thetanullwerte $\chi_{68}$,
where
$$
\chi_k(Z):=\prod_{\delta\hbox{ even}} \theta[\delta](0,Z) \ ,
$$
$Z\in\Hh_g$, with $k=2^{g-2}(2^g+1)$. Let $C$ be a Riemann surface
of genus $g\ge 2$ with a given symplectic basis for $H_1(C,\ZZ)$.
For each positive integer $n$, consider the basis
$\tilde\omega_1^{(n)},\ldots,\tilde\omega_{M_n}^{(n)}$ of $\Sym^n
H^0(K_C)$ whose elements, as in \eqref{symm}, are symmetrized tensor
products of $n$-tuples of vectors of the basis
$\omega_1,\ldots,\omega_g$, taken with respect to an arbitrary
ordering chosen once and for all. Denote by $\omega_i^{(n)}$,
$i=1,\ldots, M_n$, the image of $\tilde\omega_i^{(n)}$ under the
natural map \be\label{lapsi}\psi:\Sym^n H^0(K_C)\to H^0(K_C^n)\ .\ee
It is well known that such a map is surjective if and only if $g=2$
or $C$ is non-hyperelliptic of genus $g>2$. Of course, when
$M_n=N_n$ the map is an isomorphism. In particular, for $n=2$, this
is the case for $g=2$ and $g=3$ non-hyperelliptic. It has been shown
in \cite{BelavinTV,MorozovDA,DHokerQP} that for $g=2$
${\kuno[\omega]^{9}\over \kenne[\omega^{(2)}]} ={1\over
\pi^{12}\chi_{5}^2(\tau)}$.
Furthermore, it has been conjectured
in \cite{BelavinTV,MorozovDA} and proved in
\cite{IchikawaTI,DHokerCE} that for $g=3$ ${\kuno[\omega]^{9}\over
\kenne[\omega^{(2)}]} ={1\over 2^6\pi^{18}\chi_{18}^{1/2}(\tau)}$.

\vskip 3pt

\begin{remark}\label{vadettoqui}Under \eqref{anchequesta} we have
$\omega_i\mapsto
\omega'_i:=\sum_{j=1}^g\omega_j(C\tau+D)^{-1}_{ji}$,
$i=1,\ldots,g$. Such a transformation property induces the
$\Gamma_g$-actions $\tilde\rho^{(n)}$ on $\Sym^n H^0(K_C)$ and
$\rho^{(n)}:=\psi\circ\tilde\rho^{(n)}$ on $H^0(K_C^n)$. Explicitly,
$$\rho^{(n)}(\gamma)\cdot(\omega_{i_1}\cdots\omega_{i_n})= \sum_{j_1,\ldots,j_n=1}^g\omega_{j_1}\cdots\omega_{j_n}(C\tau+D)^{-1}_{j_1i_1}\cdots(C\tau+D)^{-1}_{j_ni_n}\ ,
$$
$\gamma\equiv\bigl({}^A_C{}^B_D\bigr)\in \Gamma_g$,
$i_1,\ldots,i_n=1,\ldots,g$. Furthermore, by \eqref{wedsym}
\be\label{implica} {\tilde{\omega}'}{}^{(n)}_{1}
\wedge\cdots\wedge{\tilde{\omega}'}{}^{(n)}_{{M_n}}=
\det(C\tau+D)^{-{g+n-1\choose n-1}} {\tilde\omega}^{(n)}_{1}
\wedge\cdots\wedge{\tilde\omega}^{(n)}_{M_n} \ . \ee\end{remark}

\vskip 3pt

\begin{definition}\label{definizionivarie}Let $n\ge 2$ be an integer. For each $i_{1},\ldots,i_{K_n}\in\{1,\ldots,M_n\}$ and for all
$x_1,\ldots,x_{N_n}\in C$, define $[i_{1},\ldots,i_{K_n}|\tau]$ to
be completely antisymmetric in $i_{1},\ldots,i_{K_n}$ and such that,
for any permutation $\pi$ of $M_n$ objects
\begin{align}[\pi(N_n+1),&\ldots,\pi(M_n)|\tau]\\ &:=
{\sgn(\pi)\det\nolimits_{1\le i,j\le N_n}
\omega^{(n)}_{\pi(i)}(x_j)\over
\kuno[\omega]^{(2n-1)^2}\theta\bigl(\sum_{1}^{N_n}
x_j-(2n-1)\Delta\bigr)\prod_{1}^{N_n}\sigma(x_j)^{2n-1}\prod_{j<k}^{N_n}
E(x_j,x_k)}\ .\notag\end{align}\end{definition}

\vskip 3pt

\noindent$[i_{1},\ldots,i_{K_n}|\tau]$ is independent of $\pi$ and of
$x_1,\ldots,x_{N_n}\in C$; in fact, analogously to Proposition \ref{thdettheta}, one can check that it is a meromorphic function with zero divisor with respect to each $x_i$ (see \cite{FayMAM,jfayy}). This definition allows to express the
generators of the kernel of the map $\psi$ in \eqref{lapsi} in terms
of the basis $\tilde \omega_1^{(n)},\ldots,\tilde
\omega_{M_n}^{(n)}$ in a very simple form.
\vskip 3pt

\begin{proposition}\label{evenstillmore}For each integer $n\geq 2$ and for all $i_{2},\ldots,i_{K_n}\in\{1,\ldots,M_n\}$ we have
\be\label{diciamolopureAAA}\sum_{i=1}^{M_n}[i,i_{2},\ldots,i_{K_n}|\tau]\,\omega^{(n)}_{i}(x)
=0\ .\ee\end{proposition}

\vskip 3pt

\noindent {\sl Proof.} If $i_{2},\ldots,i_{K_n}$ are not pairwise
distinct, this is obvious. Otherwise, the left hand side of
\eqref{diciamolopureAAA} is proportional to
$$ \sum_{\pi}[\pi(N_n+1),\ldots,\pi(M_n)|\tau]\omega^{(n)}_{\pi(N_n+1)}(x)\ ,
$$ where the sum is over the $\pi$ in $S_{M_n}$ such that $\pi(N_n+2)=i_{2},\ldots,\pi(M_n)=i_{K_n}$.
Thus, \eqref{diciamolopureAAA} is equivalent to
$\det\nolimits_{^{i\in J}_{1\le j\le N_n+1}} \omega^{(n)}_i(x_j)=0$,
with $J=\{1,\ldots,M_n\}\setminus \{i_{2},\ldots,i_{K_n}\}$ and
$x_{N_n+1}\equiv x$. \hfill$\square$

\vskip 3pt

If we identify a non-hyperelliptic smooth $C$ with its canonical
model in $\PP H^0(K_C)$, then the relations \eqref{diciamolopureAAA}
generate the ideal of hypersurfaces of degree $n$ containing the
canonical curve.

Let $\mathcal{H}_g$ be the closure of the locus of hyperelliptic
Riemann period matrices in $\Hh_g$. Consider the basis
\be\label{hyper}\eta_j=z^{j-1}dz/w \ ,\ee $j=1,\ldots,g$ of
$H^0(K_C)$, with $C$ the hyperelliptic curve
$w^2=\prod_{j=1}^{2g+2}(z-z_j)$.

\vskip 3pt

\begin{proposition}\label{fondamentale} $[i_{N_n+1},\ldots,i_{M_n}|\tau]$ are
holomorphic not identically vanishing global sections of the bundle
$${\Bbb F}_g=(\det {\Bbb E}_g)^{c_n}\otimes(\wedge^{N_n} \Sym^n{\Bbb E}_g)^*
\cong (\det {\Bbb E}_g)^{d_n}\otimes(\wedge^{K_n} \Sym^n{\Bbb
E}_g^*) \ , $$ on $\M_g$, where
\be\label{fondaa}d_n:=6n^2-6n+1-{g+n-1\choose n-1} \ ,\ee vanishing
precisely when $\omega^{(n)}_{i_1},\ldots,\omega^{(n)}_{i_{N_n}}$ is
not a basis of $H^0(K^n_C)$. In particular, they have zeroes of
order at least $(n-1)(g-1)-1$ at $\tau\in\mathcal{H}_g$.
\end{proposition}

\vskip 3pt

\noindent{\sl Proof.} Comparing Definition \ref{definizionivarie}
with \eqref{detthetaii} and \eqref{diciamolo}, yields
\be\label{mumfun} [i_{N_n+1},\ldots,i_{M_n}|\tau]=
{\epsilon_{i_1,\ldots,i_{M_n}}\omega^{(n)}_{i_1}\wedge\cdots
\wedge\omega^{(n)}_{i_{N_n}}\over (\omega_{1}\wedge\cdots
\wedge\omega_{g})^{c_n}\mu_{g,n}}\ , \ee
$i_1,\ldots,i_{M_n}\in\{1,\ldots,M_{n}\}$, where
$\epsilon_{i_1,\ldots,i_{M_n}}$ is the completely antisymmetric
tensor with $\epsilon_{1,\ldots,M_n}=1$. Holomorphicity follows by
the fact that $\mu_{g,n}$ is holomorphic on $\M_g$, and since it is nowhere vanishing on $\M_g$, even its inverse is holomorphic on $\M_g$.
In
particular, \eqref{mumfun} shows that
$[i_{N_n+1},\ldots,i_{M_n}|\tau]$ and
$\omega^{(n)}_{i_1}\wedge\cdots \wedge\omega^{(n)}_{i_{N_n}}$ have
the same divisor, so that $[i_{N_n+1},\ldots,i_{M_n}|\tau]$ does not
vanish identically on $\M_g$. Furthermore, since
$\eta_g^n=z^{n(g-1)}(dz)^n/w^n$, it follows that on the
hyperelliptic loci
 $\Sym^n H^0(K_C)$ is a $n(g-1)+1$ dimensional subspace of $H^0(K_C^n)$, so that, by \eqref{mumfun}
$[i_{1},\ldots,i_{K_n}|\tau]$ vanishes at order at least
$(n-1)(g-1)-1$ at $\tau\in\mathcal{H}_g$.
 The modular properties of
$\epsilon_{i_1,\ldots,i_{M_n}}\omega^{(n)}_{i_1}\wedge\cdots
\wedge\omega^{(n)}_{i_{N_n}}$ are the same as
$\epsilon_{i_1,\ldots,i_{M_n}}\tilde\omega^{(n)}_{i_1}\wedge\cdots
\wedge\tilde\omega^{(n)}_{i_{N_n}}$, which, in turn, can be derived
explicitly considering the identity
$$ \tilde\omega^{(n)}_1\wedge\cdots\wedge\tilde\omega^{(n)}_{M_n}=
\sum_{i_{N+1},\ldots,i_{M_n}=1}^{M_n}\epsilon_{i_1,\ldots,i_{M_n}}(\tilde\omega^{(n)}_{i_1}
\wedge\cdots \wedge\tilde\omega^{(n)}_{i_{N_n}})\otimes
\tilde\omega^{(n)}_{i_{N_{n+1}}}\otimes\cdots\otimes\tilde\omega^{(n)}_{i_{M_n}}\
.
$$ Noting that under \eqref{anchequesta}
$$ {\tilde\omega'{}^{(n)}_1\wedge\cdots\wedge\tilde\omega'{}^{(n)}_{M_n}\over
(\omega'_{1}\wedge\cdots \wedge\omega'_{g})^{c_n}\mu_{g,n}}=\det
(C\tau+D)^{d_n}{\tilde\omega^{(n)}_1\wedge\cdots\wedge\tilde\omega^{(n)}_{M_n}\over
(\omega_{1}\wedge\cdots \wedge\omega_{g})^{c_n}\mu_{g,n}}\ ,
$$
we obtain
\begin{align*}\sum_{i_{N+1},\ldots,i_{M_n}=1}^{M_n}&[i_{N_n+1},\ldots,i_{M_n}|\gamma\cdot\tau]\, \tilde\omega'{}^{(n)}_{i_{N_{n+1}}}\otimes\cdots\otimes\tilde\omega'{}^{(n)}_{i_{M_n}}\\ =&
\det
(C\tau+D)^{d_n}\sum_{k_{N+1},\ldots,k_{M_n}=1}^{M_n}[k_{N_n+1},\ldots,k_{M_n}|\tau]\,
\tilde\omega^{(n)}_{k_{N_{n+1}}}\otimes\cdots\otimes\tilde\omega^{(n)}_{k_{M_n}}\end{align*}
where $\gamma\equiv\bigl({}^A_C{}^B_D\bigr)\in \Gamma_g$,
$\gamma\cdot\tau$ given in \eqref{modull}. It follows that, under
\eqref{anchequesta},
\begin{align}
\label{transfC}\sum_{i_{1},\ldots,i_{K_n}=1}^{M_n}\rho^{(n)}(\gamma)_{k_{1}i_{1}}\cdots\rho^{(n)}(\gamma)_{k_{K_n}i_{K_n}}&[i_{1},\ldots,i_{K_n}|\gamma\cdot\tau]\\
& = \det (C\tau+D)^{d_n} [k_{1},\ldots,k_{K_n}|\tau]\ ,\notag
\end{align} so that $[i_1,\ldots,i_{K_n}|\tau]$ defines a section of the vector bundle
$$ (\det {\Bbb E}_g)^{d_n}\otimes(\wedge^{K_n} \Sym^n{\Bbb E}_g^*)\cong(\det {\Bbb E}_g)^{c_n}\otimes(\wedge^{N_n} \Sym^n{\Bbb E}_g)^* \ .$$  Here we used the isomorphisms
$$ (\wedge^{N_n}V)^*\cong \wedge^{M_n} V^*\otimes \wedge^{K_n}V^*\ ,
$$ and
$$ \wedge^{M_n} \Sym ^n W^* \cong (\det W)^{-{g+n-1\choose n-1}}\ ,
$$ that holds for a generic complex vector bundle $V$ and $W$ of rank $M_n$ and $g$, respectively (see also Proposition \ref{endo}).
 \hfill$\square$

\section{Discriminant of the $g=4$ quadric and the Schottky-Igusa
 modular form}\label{secquadric}

Denotes by $\mathcal{I}_g$ the closure of the locus of Riemann period
matrices in $\Hh_g$. The elements of such a locus can be naturally
identified with the elements of the Torelli space $\T_g$. The
Torelli space is the quotient of the Teichm\"uller space by the
Torelli group, that is the normal subgroup of the mapping class
group whose elements act trivially on $H_1(C,\ZZ)$. Equivalently,
the Torelli group is the kernel of the natural homomorphism of the
mapping class group into the symplectic group $Sp(2g,\ZZ)$. Note
that the hyperelliptic sublocus $\mathcal{H}_g\subset \Hh_g$ can be
naturally identified with the sublocus of $\T_g$ corresponding to
hyperelliptic Riemann surfaces. We denote by $\hat\T_g$ its
complement in the Torelli space, whose elements correspond to
non-hyperelliptic surfaces.

The explicit expression of the coefficients of a quadric containing
a canonical curve of genus $g$ obviously depends on the choice of a
coordinate basis of $\PP^{g-1}$ or, equivalently, of a basis of
$H^0(K_C)$. Therefore, it is natural to look for quantities
characterizing such a curve that are invariant under the projective
linear group $\PP {\rm GL}(g,\CC)$ of coordinate changes on
$\PP^{g-1}$. We denote by $I_k$ an invariant of weight $k$, i.e. a
function of the coefficients $C_{ij}$ of the quadric, transforming
as
$$I_k(C)=\det \rho(A)^kI_k(A\cdot C)\ ,$$
where $A$ is an element of ${\rm GL}(g,\CC)$, $\rho:{\rm
GL}(g,\CC)\to \End(\CC^g)$ is the fundamental representation and
$A\cdot C$ denotes the action of $A\in{\rm GL}(g,\CC)$ on the
coefficients $C_{ij}$.

If a symplectic basis of $H_1(C,\ZZ)$ is fixed, such an invariant
can be evaluated with respect to the basis of holomorphic abelian
differentials canonically normalized with respect to the
$\alpha$-periods. It follows, by definition, that $I_k$ must
transform as
$$I_k\to I_k'=\det (C\tau+D)^kI_k\ ,$$ under a $\Gamma_g$-transformation
corresponding to a change of the symplectic basis of $H_1(C,\ZZ)$.

For $g=4$ and $n=2$ Proposition \ref{evenstillmore} gives the
quadric \be\label{laequazione}
\sum_{i,j=1}^4{1+\delta_{ij}\over2}[(ij)|\tau]\omega_i\omega_j=0 \
.\ee Here, $(ij)$ denotes the element in $\{1,\ldots,M_2(4)=10\}$
such that $\omega^{(2)}_{(ij)}=\omega_i\omega_j$ in the chosen
ordering for $\omega^{(2)}_{1},\ldots,\omega^{(2)}_{10}$.

\noindent Define
\be\label{typs}\Delta_4(\tau):=\det_{ij}\big({1+\delta_{ij}\over2}[(ij)|\tau]\big)\
.\ee
 For all $i,j=1,\ldots,4$ and $Z\in \Hh_4$, set
\be\label{SISVVMF}{S_4}_{ij}(Z):={1+\delta_{ij}\over2}{\partial
F_4(Z)\over
\partial Z_{ij}} \ ,\ee
where \be\label{igusadef}F_g(Z)=2^g \sum_{\delta\hbox{
even}}\theta^{16}[\delta](0,Z)- \bigl(\sum_{\delta\hbox{
even}}\theta^{8}[\delta](0,Z)\bigr)^2 \ , \ee $Z\in\Hh_g$, is a
modular form of weight $8$. Up to normalization, it can be shown
that $F_g$ is the difference of the theta series of the even
unimodular lattices $E_8\oplus E_8$ and $D_{16}^+$
\be\label{difference}
F_g(Z)=2^{2g}(\Theta_{D_{16}^+}(Z)-\Theta_{E_8}^2(Z))\ . \ee $F_4$
is the Schottky-Igusa form \cite{Igusauno,IgusaSc} and the
irreducible variety in $\Hh_4$ defined by $F_4=0$ is $\mathcal{I}_4\subset\Hh_4$.

For each $\tau\in \I_4$, let $\Theta_s$ be the singular locus of the
theta function, and set
$$\sigma_{ij}(e,\tau):={1+\delta_{ij}\over2}{\partial\theta(e,Z)\over\partial Z_{ij}}|_{Z=\tau} \ . $$

\vskip 3pt

\begin{proposition}\label{erdiscrante} $\Delta_4(\tau)$ is a Teichm\"uller modular form
of weight $34$.\end{proposition}

\vskip 3pt

\noindent{\sl Proof.} By \eqref{transfC}, in the case $n=2$ and
$g=4$ \be\label{quattroni}[(ij)|\gamma\cdot\tau]
 = \det (C\tau+D)^8\, {}^t(C\tau+D)[(kl)|\tau](C\tau+D)\ ,\ee
so that $\Delta_4(\gamma\cdot\tau)=\det
(C\tau+D)^{34}\Delta_4(\tau)$. \hfill$\square$

\vskip 3pt

\begin{lemma}\label{lFthsing} Let $C\in \T_4$ be a marked Riemann surface
and let $\tau$  be its period matrix. Then
\be\label{Fthsing}{1+\delta_{ij}\over2}[(ij)|\tau]=c{S_4}_{ij}(\tau)
\ ,\ee $i,j=1,\ldots,4$, with $c\in\CC^\star$ independent of $\tau$,
so that \be\label{buonaa}\Delta_4(\tau)=c^4 \det{S_4}(\tau)\ .\ee
\end{lemma}

\vskip 3pt

\noindent{\sl Proof.} Define some local coordinates $t_1,\ldots,t_9$
on $\T_4$ centered at the point $C$, corresponding to the period
matrix $\tau\in\I_4$, and consider an arbitrary element
$\partial_t\in T_C\T_4$ in the tangent space. Since $F_4$ vanishes
identically on $\mathcal{I}_4$, we have
\be\label{Fgenquadr}0=\partial_t F_4(\tau)=\sum_{i\le j}{\partial
F_4\over
\partial Z_{ij}}\big\vert_{ Z=\tau}\partial_t\tau_{ij}=\sum_{i\le j}{\partial F_4\over
\partial Z_{ij}}\big\vert_{ Z=\tau}d\tau_{ij}(\partial_t)\ .\ee Here, $d\tau_{ij}$ is the element
of the cotangent space $T_C^*\T_4$ defined by
$d\tau_{ij}(\partial_t):=\partial_t\tau_{ij}$, for all
$\partial_t\in T_C\T_4$. The Kodaira-Spencer map establishes an
isomorphism between $T_C^*\T_4$ and $H^0(K_C^2)$. In particular, via
the Rauch's variational formula \cite{FayMAM}, $d\tau_{ij}$
corresponds to the quadratic differential $\omega_i\omega_j$. Since
the identity \eqref{Fgenquadr} holds for an arbitrary $\partial_t\in
T_C\T_4$, it follows that
$$\sum_{i\le j}{\partial
F_4\over
\partial Z_{ij}}\big\vert_{ Z=\tau}\omega_i\omega_j=\sum_{i,j=1}^4{S_4}_{ij}(\tau)\omega_i\omega_j=0\ ,
$$
as an element in $H^0(K_C^2)$. Since the ideal of quadrics of a
canonical curve is generated by Eq.\eqref{laequazione}, we have
${S_4}_{ij}(\tau)=c(\tau){1+\delta_{ij}\over2}[(ij)|\tau]$, for some
holomorphic function $c(\tau)$ on $\mathcal{I}_4$, independent of
$i,j=1,\ldots,4$. Let us prove that $c(\tau)$ must be invariant
under the action of $\Gamma_4$ on $\tau$. Since $F_4(\tau)=0$ for
all $\tau\in \I_4$, it follows that on $\I_4$
\be\label{Fijtransf}{S_4}(\gamma\cdot\tau)= \det(C\tau+D)^8\,
{}^t(C\tau+D) {S_4}(\tau)(C\tau+D)\ ,\ee which is the same
transformation property satisfied by $[(ij)|\tau]$. Thus, $c$ is
modular invariant, so it must be a constant. Finally, observe that
$c$ cannot vanish since it would imply ${\partial F_4\over\partial
Z_{ij}}(\tau)=0$ for all $\tau\in\mathcal{I}_4$, which is impossible
because $F_4$ is irreducible \cite{IgusaSc}.\hfill$\square$

\vskip 3pt

\noindent It is now useful to recall the following well known
result.

\vskip 3pt

\begin{lemma}\label{lethetanull}Let $C$ be either a non-hyperelliptic Riemann surface of genus $g=4$
or a non-trigonal surface of $g=5$. Then, the canonical model of $C$
is contained in a quadric of rank $3$ if and only if
$\prod_{\delta\hbox{ even}}\theta[\delta]=0$.\end{lemma}

\vskip 3pt

\noindent{\sl Proof.} The modular form $\prod_{\delta\hbox{
even}}\theta[\delta]$ vanishes if and only if $C$ has an even
singular spin structure $\delta$. In this case, there are two
holomorphic sections $\xi_1$ and $\xi_2$ of $L_\delta$,
$L_\delta^2=K_C$, so that taking $\eta_1=\xi_1^2$, $\eta_2=\xi_2^2$
and $\eta_3=\xi_1\xi_2$ one has $\eta_3^2 = \eta_1\eta_2$ which is a
quadric of rank $3$ containing $C$.

Conversely, suppose that $\eta_3^2 = \eta_1\eta_2$ for some
$\eta_1,\eta_2,\eta_3\in H^0(K_C)$. Set $(\eta_i)=\sum_{p\in
C}m_i(p)p$ and consider the divisor $D=\sum_{p\in
C}\min\{m_1(p),m_3(p)\}p$. $D$ has degree at most $g-1$, otherwise
the ratio $\eta_1/\eta_3$ would be a meromorphic function with at
most $g-2$ poles and the curve would be hyperelliptic for $g=4$ or
trigonal for $g=5$. On the other hand, $\eta_3^2 = \eta_1\eta_2$
implies that $(\eta_1)\leq 2(\eta_3)$ and since the supports of
$(\eta_1)-D$ and $(\eta_3)-D$ are disjoint, the only possibility is
that $(\eta_1)=2D$. Therefore, $\eta_1$ is the square of a
holomorphic section of the line bundle $L_\delta$, with
$L_\delta^2=K_C$, corresponding to the divisor $D$. By the same
reasoning it follows that $\eta_2$ is the square of a holomorphic
section of the line bundle $L_{\delta'}$, with $L_{\delta'}^2=K_C$,
corresponding to the divisor $D'$. Since $(\eta_3)=D+D'$ is a
canonical divisor it follows that $\delta'=\delta$. Then $\delta$
necessarily corresponds to an even singular spin structure, since a
surface of genus $g=4,5$ admits no odd super-singular spin
structures, that is spin structures with three or more holomorphic
sections.\hfill$\square$

\vskip 3pt

\begin{theorem}\label{thdettau} For any $\tau\in\mathcal{I}_4$
\be\label{sqdet}\det S_4(\tau)=d\,\chi_{68}(\tau)^{1/2}\ ,\ee  where
$d\in\CC^\star$ is independent of $\tau$. In particular, by Lemma
\ref{lFthsing}, \be\label{vspritta}\Delta_4(\tau)={d\over
c^4}\chi_{68}(\tau)^{1/2} \ .\ee
Furthermore,
\be\label{numbertheory}
\chi^{68}(\tau)^{1/2}={2^{28}\over d}\det[(1+\delta_{ij}){\partial\over\partial
Z_{ij}}(\Theta_{D_{16}^+}-\Theta_{E_8}^2)|_{Z=\tau}] \ .
\ee
\end{theorem}

\vskip 3pt

\noindent {\sl Proof.} By \eqref{Fijtransf} it follows that $\det
S_4(\gamma\cdot\tau)= \det(C\tau+D)^{34}\det S_4(\tau)$ so that
$\det S_4$ is a modular form of weight $34$ when restricted to
$\mathcal{I}_4$. (Note that $\det S_4$ is not a modular form on the
whole $\Hh_4$, since the modular group action on $S_4$ is affine
outside $\I_4$). On the other hand, it has been proved in
\cite{tsuy} that the square root of $\chi_{68}$ in the RHS of
\eqref{sqdet} is well-defined when restricted to $\mathcal{I}_4$.

Choose $\tau\in\mathcal{I}_4\setminus \mathcal{H}_4$. By Lemma
\ref{lFthsing}, the LHS of \eqref{sqdet} is proportional to the
discriminant. By Lemma \ref{lethetanull}, $\Delta_4(\tau)$ vanishes
if and only if the Riemann surface has a singular even spin
structure. The locus of Riemann surfaces with singular even spin
structures in $\T_4$ corresponds to the divisor of
$\sqrt{\prod_{\delta\hbox{ even}}\theta[\delta](0,\tau)}$ in $\mathcal{
I}_4$ \cite{tsuy}, so that the meromorphic function
$${\det{ S_4}\over \sqrt{\prod_{\delta\hbox{
even}}\theta[\delta](0,\tau)}}\ ,$$ has no poles on $\mathcal{
I}_4\setminus \mathcal{H}_4$ and therefore on $\mathcal{I}_4$, since
$\mathcal{H}_4$ has codimension $2$. Since it is a holomorphic modular
invariant function not identically zero on $\mathcal{I}_4$, it must be
a non-vanishing constant. Finally, \eqref{numbertheory} follows from (\ref{difference}) and (\ref{sqdet}). \hfill$\square$

\vskip 3pt

\begin{theorem}\label{thetasingularfour} Let $e\in\Theta_s$. The following functional relation between
$\Theta_s$ and $\tau\in \mathcal{I}_4$ \be\label{thetasingeF}
S_{4ij}(\tau)=\Bigg({d\chi_{68}(\tau)^{1/2}\over\det\sigma(e,\tau)}\Bigg)^{1/4}\sigma_{ij}(e,\tau)
\ , \ee  holds. Both $S_{4ij}(\tau)$ and $\sigma_{ij}(e,\tau)$ have
a zero of order at least $2$ for all $\tau\in\mathcal{H}_4$ for any
$i,j=1,\ldots,4$. Furthermore, both $\det S_4$ and
$\det\sigma(e,\tau)$ have a zero of order 10 for all $\tau\in\mathcal{
H}_4$ and at least a zero when $e\in\Theta_s$ is an even
$\theta$-characteristic for all $\tau\in\I_4$.

\end{theorem}

\vskip 3pt

\noindent {\sl Proof.}  The well known relation
$\sum_{i,j=1}^4\sigma_{ij}(e,\tau)\omega_i\omega_j=0$,
$e\in\Theta_s$, implies
\be\label{bbbg}S_{4ij}(\tau)=G(e|\tau)\sigma_{ij}(e,\tau) \ ,\ee
with $G(e|\tau)$ determined by Theorem \ref{thdettau} upon taking
the determinant of both sides of \eqref{bbbg}. The fact that all the
$\tau$ in $\mathcal{H}_4$ are zeros of order at least $(2-1)(4-1)-1=2$
of $S_{4ij}(\tau)$ for any $i,j=1,\ldots,4$, is an immediate
consequence of Proposition \ref{fondamentale} and  Lemma
\ref{lFthsing}. Concerning the order of the zero of $\det S_4$ and
$\det\sigma(e,\tau)$ for all $\tau\in\mathcal{H}_4$, it immediately
follows by (\ref{thetasingeF}) and the well known fact that for each
$\tau\in\mathcal{H}_4$ there are 10 vanishing thetanulls. Finally,
since the order of the zeros of $S_{4ij}(\tau)$ depends on $i,j$ and
$\tau$, it follows by \eqref{thetasingeF} that $S_{4ij}(\tau)$ and
$\sigma_{ij}(e,\tau)$ have the same divisors in $\I_4$. When
$e\in\Theta_s$ is an even characteristic $\chi_{68}(\tau)$ has at
least a double zero, so that $\det S_4(\tau)$ and by
\eqref{thetasingeF}, $\det\sigma(e,\tau)$, has a zero of order at
least one for all $\tau\in\I_4$.  \hfill$\square$

\vskip 3pt

\begin{remark}\label{vadettoanchequa} In \cite{LRZ} it has been re-obtained the Klein formula linking
$\chi_{18}$ to the square of the discriminant of plane quartics. The
authors also studied possible generalizations to the case $g>3$. In
particular, in Eq.(3) of \cite{LRZ}, it has been mentioned the
following ``amazing formula'' by Klein in the footnote, p.462 in
\cite{KLEIN} \be\label{Klein} \chi_{68}(\tau)=c'\tilde\Delta_4(C)^2
T(C)^8 \ , \ee with $c'$ a constant. This formula relates
$\chi_{68}$ to the discriminant $\tilde\Delta_4(C)$ of the Klein
quadric and the tact invariant $T(C)$
 of the quadric and of the cubic (see pg.122 of \cite{Salomon}), whose intersection in $\PP^3$ determines $C$. Note that
if $\sum_{i,j=1}^4C_{ij}\omega_i\omega_j=0$, denotes  the Klein
quadric, then \be\label{relaK}C_{ij}=\tilde c {S_{4ij}\over T(C)} \
,\ee with $\tilde c\in\CC^\star$ independent of $\tau$. \end{remark}

\vskip 3pt

The following expression for $\mu_{4,2}$ has been suggested in
\cite{BelavinGA} and in \cite{MorozovDA} in the context of bosonic
string theory. Its proof has been a long-standing problem and a more
rigorous derivation has been provided in \cite{Guilarte}. In the
present approach it follows immediately.

\vskip 3pt

\begin{theorem}\label{stillmore}The $g=4$, $n=2$ Mumford form is
\be\label{mumffffour}\mu_{4,2}=\pm{1\over c
S_{4ij}}{\omega_1\omega_1\wedge\cdots\wedge
\widehat{\omega_i\omega_j}\wedge\cdots\wedge \omega_4\omega_4\over
(\omega_1\wedge\cdots\wedge\omega_4)^{13}} \ ,\ee for all
$i,j=1,\ldots,4$, with $c\in\CC^*$ the constant defined in
Eq.\eqref{Fthsing}.\end{theorem}

\vskip 3pt

\noindent {\sl Proof.} Immediate by  Eq.\eqref{mumfun} and
Eq.\eqref{Fthsing} of Lemma \ref{lFthsing}. \hfill$\square$

\newpage

\bibliographystyle{amsplain}

\begin{thebibliography}{99}

\bibitem{Matone:2005vm}
  M.~Matone and R.~Volpato,
  Higher genus superstring amplitudes from the geometry of moduli space,
  Nucl.\ Phys.\ B {\bf 732} (2006) 321
  [hep-th/0506231].

\bibitem{Matone:2008td}
  M.~Matone and R.~Volpato,
  Superstring measure and non-renormalization of the three-point amplitude,
  Nucl.\ Phys.\ B {\bf 806} (2009) 735
  [arXiv:0806.4370 [hep-th]].

\bibitem{Matone:2010yv}
  M.~Matone and R.~Volpato,
  Getting superstring amplitudes by degenerating Riemann surfaces,
  Nucl.\ Phys.\ B {\bf 839} (2010) 21
  [arXiv:1003.3452 [hep-th]].

\bibitem{Matone:2005bx}
  M.~Matone and R.~Volpato,
  Linear relations among holomorphic quadratic differentials and induced Siegel's metric on M(g),
  J.\ Math.\ Phys.\  {\bf 52} (2011) 102305
  [math/0506550 [math.AG]].

\bibitem{Matone:2012wy}
  M.~Matone,
  Extending the Belavin-Knizhnik ``wonderful formula'' by the characterization of the Jacobian,
  JHEP {\bf 1210} (2012) 175
  [arXiv:1208.5994 [hep-th]].

\bibitem{FayMAM}
J.~Fay, Kernel functions, analytic torsion and moduli spaces, {\it
Mem.\ Am.\ Math.\ Soc.\ } {\bf 96} (1992).

\bibitem{jfayy}
J.~Fay, {\it Theta Functions on Riemann surfaces}, Springer Lecture
Notes in Math. {\bf 352}, 1973.

\bibitem{MatoneBB}
  M.~Matone and R.~Volpato,
  Determinantal characterization of canonical curves and combinatorial theta
  identities. To be published in Math. Ann.

\bibitem{Mumford}
D.~Mumford, Stability of projective varieties, {\it Enseign.\ Math.\
} {\bf 23} (1977), 39-110.

\bibitem{BelavinCY}
  A.~A.~Belavin and V.~G.~Knizhnik,
  Algebraic geometry and the geometry of quantum strings,
  {\it Phys.\ Lett.\ } B {\bf 168} (1986), 201-206.

\bibitem{BonoraCJ}
  L.~Bonora, A.~Lugo, M.~Matone and J.~Russo,
  A Global Operator Formalism on Higher Genus Riemann Surfaces: b-c
  Systems,
 {\it Commun.\ Math.\ Phys.\ } {\bf 123}  (1989), 329-352.

\bibitem{BostUY}
  J.~B.~Bost and T.~Jolicoeur, A holomorphy property and critical dimension in string theory and index theorem,
   {\it Phys.\ Lett.\ } B {\bf 174} (1986), 273-276.

\bibitem{BismutWR}
J.~M.~Bismut and D.~S.~Freed, The analysis of elliptic families. I.
Metrics and connections on determinant bundles, {\it Commun.\ Math.\
Phys.\ } {\bf 106} (1986), 159-176; The analysis of elliptic
families. II. Dirac operators, eta invariants, and the holonomy
theorem,
 {\bf 107} (1986), 103-163.

\bibitem{AlvarezGaumeVM}
  L.~Alvarez-Gaume, J.~B.~Bost, G.~W.~Moore, P.~C.~Nelson and C.~Vafa,
  Bosonization on higher genus Riemann surfaces,
 {\it Commun.\ Math.\ Phys.\ } {\bf 112} (1987), 503-552.

\bibitem{RaySinger}
D.~Ray and I.~M.~Singer, Analytic torsion for complex manifolds,
{\it Ann.\ Math.\ } {\bf 98} (1973), 154-180.

\bibitem{BeilinsonZW}
  A.~A.~Beilinson and Y.~I.~Manin,
  The Mumford form and the Polyakov measure in string theory,
  {\it Commun.\ Math.\ Phys.\ }  {\bf 107} (1986), 359-376.

\bibitem{VerlindeKW}
  E.~P.~Verlinde and H.~L.~Verlinde,
  Chiral bosonization, determinants and the string partition function,
{\it  Nucl.\ Phys.\ } B {\bf 288} (1987), 357-396.

\bibitem{vanderGeer}
G.~van~der~Geer, Siegel modular forms and their applications, The
1-2-3 of modular forms, 181-245, Universitext, Springer, Berlin,
2008, [arXiv:math.AG/0605346].

\bibitem{MatoneZZ}
  M.~Matone and R.~Volpato,
  The singular locus of the theta divisor and quadrics through a canonical
  curve,
  arXiv:0710.2124 [math.AG].

\bibitem{ShepherdBarronHC}
  N.~I.~Shepherd-Barron,
  Thomae's formulae for non-hyperelliptic curves and spinorial square roots
  of theta-constants on the moduli space of curves,
  arXiv:0802.3014 [math.AG].

\bibitem{BelavinTV}
  A.~A.~Belavin, V.~Knizhnik, A.~Morozov and A.~Perelomov,
  Two and three loop amplitudes in the bosonic string theory,
  {\it Phys.\ Lett.\ } B {\bf 177} (1986), 324-328. 

\bibitem{MorozovDA}
  A.~Morozov,
  Explicit formulae for one, two, three and four loop string amplitudes,
 {\it  Phys.\ Lett.\ } B {\bf 184} (1987), 171-176.

\bibitem{DHokerQP}
E.~D'Hoker and D.~H.~Phong, Two-loop superstrings. IV: The
cosmological constant and modular forms, {\it Nucl.\ Phys.\ } B {\bf
639} (2002), 129-181, [arXiv:hep-th/0111040].

\bibitem{IchikawaTI}T.~Ichikawa, On Teichm\"uller modular forms,
{\it  Math.\ Ann.\ }  {\bf 299}  (1994), 731-740.

\bibitem{IchikawaT}
T.~Ichikawa, Teichm\"uller modular forms of degree 3, {\it Amer.\
J.\ Math.\ } {\bf 117} (1995), 1057-1061.

\bibitem{DHokerCE}
  E.~D'Hoker and D.~H.~Phong,
  Asyzygies, modular forms, and the superstring measure. II,
  {\it Nucl.\ Phys.\ } B {\bf 710} (2005), 83-116,
  [arXiv:hep-th/0411182].

\bibitem{Igusauno}
J.-I.~Igusa, Schottky's invariant and quadratic forms, {\it E. B.
Christoffel} (Aachen/Monschau, 1979), 352-362, Birkh\"auser,
Basel-Boston, Mass., 1981.

\bibitem{IgusaSc}J.-I.~Igusa, On the irreducibility of Schottky's
divisor, {\it J. Fac. Sci. Univ. Tokyo Sect. IA Math.} {\bf 28}
(1981), 531-545.

\bibitem{tsuy}
S.~Tsuyumine, Thetanullwerte on a moduli space of curves and
hyperelliptic loci, {\it Math. Z.} {\bf 207} (4) (1991), 539-568.

\bibitem{LRZ}
G.~Lachaud, C.~Ritzenthaler and A.~Zykin, Jacobians among abelian
threefolds: a formula of Klein and a question of Serre, {\it Math.\
Res.\ Lett.\ } {\bf 17} (2010), no. 2, 323-333,
[arXiv:0802.4017[math.NT]].

\bibitem{KLEIN}
F.~Klein, Zur theorie der abelschen funktionen, {\it Math.\ Annalen
} {\bf 36} (1889-90); = Gesammelte mathematische Abhandlungen XCVII,
388-474.

\bibitem{Salomon}
G.~Salmon, {\it Trait\'e de g\'eom\'etrie analytique \`a trois
dimensions}. Troisi\`eme partie. Ouvrage traduit de l'anglais sur la
quatri\`eme \'edition, Paris, 1892.

\bibitem{BelavinGA}
  A.~A.~Belavin and V.~G.~Knizhnik,
  Complex geometry and the theory of quantum strings,
  {\it Sov.\ Phys.\ JETP } {\bf 64} (1986), 214-228
  [{\it Zh.\ Eksp.\ Teor.\ Fiz.\ }  {\bf 91} (1986), 364-390].

\bibitem{Guilarte}
J.~M.~Guilarte and J.~Mu{\~n}oz Porras, Four-loop vacuum amplitudes
for the bosonic string, {\it Proc. \ Roy. \ Soc. \ London \ } Ser. A
{\bf 451} (1995), 319-329.

\end{thebibliography}

\end{document}